\documentclass[11pt,reqno]{amsart}

\usepackage[margin=1in]{geometry}

\usepackage{amsmath}
\usepackage{amsfonts}
\usepackage{amssymb}
\usepackage{bm} 

\newcommand\R{{\mathbb{R}}}

\renewcommand{\div}{\operatorname{div}}

\theoremstyle{plain}
  \newtheorem{theorem}[subsection]{Theorem}

  \newtheorem{corollary}[subsection]{Corollary}

\theoremstyle{remark}
  \newtheorem{remark}[subsection]{Remark}
  
    \newtheorem{example}[subsection]{Example}	

\theoremstyle{definition}

\include{psfig}

\begin{document}


\title[Uniqueness of weak solutions to active scalar equations]{Logarithmic spikes of gradients and uniqueness of weak solutions  to a class of active scalar equations}

\author{Walter Rusin}
\address{Department of Mathematics \\ University of Southern California \\ 3620 S. Vermont Ave. \\ KAP 108 \\ Los Angeles, CA 90089}
\email{wrusin@usc.edu}
\date{\today}

\begin{abstract} 
	We study the question weather weak solutions to a class of active scalar equations, with the drift velocity and the active scalar related via a Fourier multiplier of order zero, are unique. Due to some recent results we cannot expect weak solutions to be unique without additional conditions. We analyze the case of some integrability conditions on the gradient of the solutions. The condition is weaker than simply imposing $\nabla \theta \in L^\infty$. Lastly, we consider the inviscid limit for the studied class of equations. 
\end{abstract}

\maketitle

\section{Introduction}
Active scalar evolution equations have been a topic of considerable study in recent years, in part because they arise in many physical models.  In particular, such equations are prevalent in fluid dynamics. This paper is concerned with a class of diffusive transport equations for an unknown scalar field $\theta(x,t)$, of the form
\begin{equation}\label{DEQ}\tag{DDE}
	\partial_t \theta + (u \cdot \nabla) \theta  + \kappa(-\Delta)^\gamma \theta= 0,
\end{equation}
where $\kappa>0$, and their non-diffusive version
\begin{equation}\label{EQ}\tag{DE}
	\partial_t \theta + (u \cdot \nabla) \theta  = 0,
\end{equation}
For both equations we assume $t>0$ and $x \in \R^n$. The divergence free velocity field $u(x,t)$
\begin{equation}
	\div u = 0
\end{equation}
is obtained from $\theta$ via Calder\'on-Zygmund singular integral operators, defined in the frequency space by a multiplier, whose symbol is of order zero (see \cite{Stein} for a precise definition)
\begin{equation}
	\hat{u}(\xi) = T(\xi)\hat{\theta}(\xi).
\end{equation}
Examples covered by the above class of equations arise in a variety of physical models. To be more specific we present a few examples of such.

\begin{example}
	The 2D surface quasi-geostrophic equation models the evolution of buoyancy (or potential temperature) $\theta$ on the 2D horizontal boundaries of general 3D quasi-geostrophic equations (see \cite{CMT94} for more details). It is given by
	$$
		\partial_t \theta + (u \cdot \nabla) \theta + (-\Delta)^\alpha \theta =0,
	$$
	with $\div u=0$. The drift velocity $u$ is recovered from the buoyancy $\theta$ through the multiplier operator
	$$
		|\xi|^{-1}(-\xi_2,\xi_1).
	$$
\end{example}

\begin{example}
	The 2D porous media equation governs the motion of the density $\rho$ of an incompressible fluid:
	$$
		\partial_t \rho + (u \cdot \nabla) \rho =0,
	$$
	with $\div u=0$, where $u$ and $\rho$ are related via Darcy's law:
	$$
		u = -\nabla p - (0,\rho).
	$$
	For simplicity we have set all physical constants equal to $1$. Eliminating the pressure from the equation, one obtains $u = T[\rho]$, where $T[\,\cdot\,]$ is a Fourier multiplier operator with symbol
	$$
		 |\xi|^{-2}(\xi_1\xi_2,-(\xi_1)^2).
	$$
\end{example}

In the present paper, we shall be concerned with \emph{weak solutions}. We say that a pair $(\theta,u) \in L^2_{loc}(\R^n \times \R)$ is a weak solution of the equation (\ref{DEQ}) on $\R$ if for every $\phi \in C^\infty_0(\R^n \times \R)$ the following holds
$$
	\int_{\R}\int_{\R^n} \theta(\partial_t \phi + (u \cdot \nabla) \phi + \kappa(-\Delta)^\gamma \phi) \;dxdt = 0.
$$
The primary motivation for studying weak solutions in the context of hydrodynamics comes from the desire to understand the laws of turbulence in the limit of infinite Reynolds number. For the above mentioned examples of active scalar equations, existence of weak solutions has been shown for instance in \cite{Res} and \cite{CCG}, respectively. For the general dissipative case addressed in this paper, we give a proof of existence of weak solutions in the appendix. Proofs of existence for the non-dissipative equations take into account the special structure of the non-linearity (the actual symbols of the multipliers and related to them commutator estimates) thus we refer the interested reader to papers devoted to the analysis of a specific problem. 

The classical way to prove uniqueness of such solutions relays on applying energy methods to an equation governing the difference of two solutions, $\tilde{\theta} = \theta_1-\theta_2$. The most direct condition which guarantees the uniqueness of $\theta$ is to require, for instance, that $\nabla \theta_1 \in L^1(0,T,L^\infty(\R^n))$. We then obtain the following immediate estimate
\begin{equation}
	\left| \int_0^T \int_{\R^n} \tilde{\theta} \cdot \nabla \theta_1  \tilde{\theta} \;dxdt \right| \leq C \|\nabla \theta_1\|_{L^1((0,T),L^\infty(\R^n))}\|\tilde{\theta}\|^2_{L^\infty((0,T),L^2(\R^n))},
\end{equation}
which is the core of the energy approach. Given this estimate, uniqueness of solutions follows from an application of Gronwall's inequality. Note that we impose high regularity on the solution (and it is not clear that it follows from energy methods). Moreover, for equations in 2D, the classical result gives local existence of solutions $\theta$ with data in Sobolev space $H^s$ with $s > 2$. Even for solutions in $H^2$, we have $\nabla \theta \in H^1$. This, however, is not enough to conclude that $\nabla \theta \in L^\infty$ and motivates the choice of $BMO$ (see next section for the definition) rather than $L^\infty$.

Results of C\'ordoba, Faraco, Gancedo \cite{CFG} and Shvydkoy \cite{S11} imply that it is not enough to require only that $(\theta,u) \in L^\infty(\R^n\times \R)$ (obviously contained in $L^2_{loc}(\R^n\times \R)$) in order to obtain unique weak solutions. Therefore, one needs additional conditions in order to ensure uniqueness of such solutions. 

In our paper we consider a situation where $\nabla \theta = g+b$ with $g \in L^1((0,T),L^\infty(\R^n))$ and $b \in L^p((0,T),BMO)$ such that $supp \; b$ has finite measure for a.e. $t \in (0,T)$.  More precisely, in Section 3, we prove that under such condition, where $p$ is related to the strength of dissipation $\gamma$, weak solutions to (\ref{DEQ}) are unique. We prove an analogous result in the case of the inviscid problem (\ref{EQ}), however here we need to impose additional integrability of solutions (since it is no longer provided by dissipation). In Section 4, we consider the inviscid limit of (\ref{DEQ}). Our proofs follow techniques introduced in \cite{MR08}.

\section{Preliminaries}

We begin by recalling the definition of the space $BMO$. Given a function $f \in L^1_{loc}(\R^n)$ and a cube $Q$, let $f_Q$ denote the average of $f$ on $Q$
$$
	f_Q = \frac{1}{|Q|} \int_Q f(x)\;dx. 
$$
Define the sharp maximal function by
$$
	M^\#f(x) = \sup_{Q} \frac{1}{|Q|} \int_Q |f(x) - f_Q|\;dx,
$$
where the supremum is taken over all cubes $Q$ containing $x$. We say that $f$ has bounded mean oscillation if the function $M^\#f$ is bounded, that is
$$
	BMO = \{ f \in L^1_{loc}(\R^n): M^\# f \in L^\infty(\R^n)\}.
$$
We define a norm on $BMO$ by
$$
	\|f\|_{BMO} = \|M^\# f\|_{L^\infty}.
$$
Clearly, $L^\infty \subseteq BMO$. Recall however that $L^\infty \subsetneq BMO$ due to the fact that BMO contains also, for instance, unbounded functions with logarithmic singularities, like
$$
	f(x) = \left\{\begin{array}{cll}
		-\ln|x| & & \text{for $0<|x|<1$} \\
		& & \\
		0 & & \text{ for $|x| \geq 1$}.
	\end{array}\right.
$$
As shown by our results, the presence of such logarithmic (or in general $BMO$) singularities, localized to a set of finite measure, does not affect the uniqueness of weak solutions.

In the proofs we relay on the following remark.

\begin{remark}\label{rem1}
	Observe that $f \in BMO$ is by definition equivalent to $M^\# f \in L^\infty$. Because of the above comment, the reverse implication is not true. However, Fefferman and Stein in~\cite{FS72} proved that $M^\# f \in L^p$ implies $f \in L^p$, if $1 < p < \infty$. Moreover, the John-Nirenberg inequality yields that if $f$ is in BMO, then $f$ is locally in $L^p$ for any $1<p<\infty$. 
\end{remark}

It is not obvious that one can consider a dual pairing of $L^1$-functions with functions in $BMO$ because $(BMO)^*=\mathcal{H}^1 \subsetneq L^1$ (the Hardy space). This problem was considered by P.B. Mucha and the author in \cite{MR08}. We recall the result since it will be the main tool used in the proofs.
\begin{theorem}\label{log_inequality}
	Let $f \in BMO$, such that $supp f$ has finite measure and $g \in L^1(\R^n)\cap~L^\infty(\R^n)$. Then the following inequality is true
	\begin{equation}
		\left| \int_{\R^n} fg \;dx \right| \leq C \|f\|_{BMO} \|g\|_{L^1} \left[1+ |\ln \|g\|_{L^1}| + \ln(1+\|g\|_{L^\infty})\right]. 
	\end{equation}
\end{theorem}

\section{Uniqueness of weak solutions}

We first state and prove the result for diffusive equations.

\begin{theorem}\label{thm1}
Let $\theta(x,t)$ be a weak solution of the equation (\ref{DEQ}) such that 
	$$
		\theta \in L^\infty((0,T),L^2(\R^n))\cap L^2((0,T),H^\gamma(\R^n)).
	$$
	Assume in addition $\nabla \theta =  g+b$ where $g \in L^1((0,T),L^\infty(\R^n))$ and $b \in L^s((0,T),BMO)$ for some $s>\frac{q}{q-2}$ where $q \geq 2$ ($s=\infty$ if $q=2$) is such that $\theta \in L^q((0,T),L^p(\R^n))$ for some $p \in (2,\frac{2n}{n-2\gamma})$, with $supp \; b$ of finite measure for a.e. $t \in (0,T)$. Then $(\theta,u)$ is unique.
\end{theorem}
\begin{proof} Let $\theta_1$ and $\theta_2$ be two solutions of (\ref{DEQ}) and $u_1, u_2$ corresponding drift velocities. The difference $\theta=\theta_1-\theta_2$ satisfies the following equation
\begin{equation}
	\theta_t + (u_1-u_2) \cdot \nabla \theta_1 + u_2\cdot \nabla \theta + \nu (-\Delta)^\gamma \theta = 0.
\end{equation}
Let us denote $\Lambda = (-\Delta)^{1/2}$. Multiplying by $\theta$ and integrating in space, using the divergence-free condition on $u_2$, we obtain
\begin{equation}\label{7}
	\frac{1}{2}\frac{d}{dt} \int_{\R^n} |\theta|^2\;dx + \nu \int_{\R^n} |\Lambda^{\gamma}\theta|^2\;dx + \int_{\R^n}(u_1-u_2)\cdot\nabla \theta_1 \theta\;dx =0.
\end{equation} 	
Set $\alpha=(u_1-u_2)\theta$ and $\beta=\nabla \theta_1$. Notice that $\|\theta\|_{L^2}$ and $\|u\|_{L^2}$ are both $C([0,T])$ upon (\ref{7}) and assumptions on $\theta_1$ and $\theta_2$ . We split $\alpha = \alpha_m + \alpha_r$, where $|\alpha_m| = \min(|\alpha|,m)$ for some $m>1$. Let $\beta = \beta_1 + \beta_2$, where $\beta_1$ corresponds to the part $g$ and $\beta_2$ corresponds to the part $b$ of the gradient.Theorem \ref{log_inequality} yields
\begin{equation}
	\begin{gathered}
	\left| \int_{\R^n} |\alpha\beta|\;dx\right| \hfill\\
	 \leq \|\beta_1\|_{L^\infty}\|\alpha\|_{L^1}+C\|\beta_2\|_{BMO}\|\alpha_m\|_{L^1} \left( 1+ |\ln\|\alpha_m\|_{L^1}|+ \ln(1+m)\right) + \int_{\R^n}|\alpha_r\beta_2|\;dx.
	 \end{gathered}
\end{equation}
Notice that $\|\alpha_m\|_{L^1} \leq \|\alpha\|_{L^1} \leq \|u\|_{L^2}\|\theta\|_{L^2}\leq C_1\|\theta\|^2_{L^2}$. Because of continuity of $\|\theta\|_{L^2}$ and $\|u\|_{L^2}$ in time, we can pick a time $T_0>0$ such that for $t \in [0,T_0]$ 
\begin{equation}
	\|\alpha_m\|_{L^1}|\ln\|\alpha_m\|_{L^1}| \leq C_1\|\theta\|^2_{L^2}|\ln C_1\|\theta\|^2_{L^2}|,
\end{equation}
thus we get
\begin{equation}\label{10}
	\begin{gathered}
	\left| \int_{\R^n} |\alpha\beta|\;dx\right| \hfill\\
	\leq C(\|\beta_1\|_{L^\infty}+\|\beta_2\|_{BMO})C_1\|\theta\|^2_{L^2} \left( 1+ |\ln C_1\|\theta\|^2_{L^2}|+ \ln(1+m)\right) + \int_{\R^n}|\alpha_r\beta_2|\;dx.
	\end{gathered}
\end{equation}
Inequality (\ref{10}) combined with (\ref{7}), after dropping the positive dissipation term, give
\begin{equation}
	\begin{gathered}
	\frac{1}{2C_1}\frac{d}{dt}C_1\|\theta\|^2_{L^2} \hfill \\
	\leq C(\|\beta_1\|_{L^\infty}+\|\beta_2\|_{BMO})C_1\|\theta\|^2_{L^2} \left( 1+ |\ln C_1\|\theta\|^2_{L^2}|+ \ln(1+m)\right) + \int_{\R^n}|\alpha_r\beta_2|\;dx.
	\end{gathered}
\end{equation}
Denote
$$
	f(t) = 2C_1C(\|\beta_1\|_{L^\infty}+\|\beta_2\|_{BMO}), \;\;\;\;\; r(t)=C_1\int_{\R^n}|\alpha_r\beta_2|\;dx, \;\;\;\;\; x(t) = C_1\|\theta\|^2_{L^2}.
$$
Notice that the assumptions on $\nabla \theta$ guarantee that $f(t)$ is integrable on bounded time intervals. We obtain the following differential inequality
\begin{equation}
	\begin{gathered}
	x'(t) \leq f(t) x(t) (1+|\ln x(t)| + \ln(1+m))+r(t),  \\
	x(0)=0. \hfill
	\end{gathered}
\end{equation}
In order to find a good estimate on $x(t)$ we introduce the following equation
\begin{equation}\label{2.10}
	\begin{gathered}
	y'(t) = f(t) y(t) (1+|\ln y(t)| + \ln(1+m))+r(t),  \\
	y(0)=1/m, \hfill
	\end{gathered}
\end{equation}
for some $m$ large enough. From the Osgood existence theorem we know that there exists a unique local solution to (\ref{2.10}). In addition, the right hand-side of (\ref{2.10}) guarantees that $y(t)$ is increasing. This implies that the solution of (\ref{2.10}) dominates $x(t)$, that is 
$$
	0 \leq x(t) \leq y(t) \leq 1 \;\;\;\;\; \text{ for $t \in [0,T_0]$.}
$$  
Therefore, we investigate the behavior of solutions to (\ref{2.10}). By Gronwall's inequality we obtain
\begin{equation}\label{14}
	\begin{gathered}
		y(t) \leq \frac{1}{m} \exp\left( \int_0^t f(s)(1+|\ln y(s)|+\ln(1+m))\;ds \right)   \\
		+ \int_0^t r(s) \exp\left( \int_s^t f(\tau)(|\ln y(\tau)| + 1 + \ln(1+m))\; d\tau  \right)\;ds \\
		\leq \frac{1}{m}\exp\left((1+\ln(1+m))\int_0^tf(s)\;ds \right)\exp\left( \int_0^tf(s)|\ln y(s)|\;ds\right) \\
		+ \exp\left( \int_0^tf(\tau)|\ln y(\tau)|\;d\tau|\right)\exp\left((1+\ln(1+m))\int_0^tf(\tau)\;d\tau \right)\int_0^tr(s)\;ds.
	\end{gathered}
\end{equation}
Since $1 \geq y(t) \geq 1/m$ implies $|\ln y(t)| \leq \ln m$ we can estimate the right hand-side of (\ref{14}) (up to a constant) by
\begin{equation}
	\begin{gathered}
		\frac{1}{m}(1+m)^{\int_0^tf(s)\;ds}m^{\int_0^tf(s)\;ds}+ m^{\int_0^tf(s)\;ds}(1+m)^{\int_0^tf(s)\;ds}\int_0^tr(s)\;ds \\
		\leq (2m^2)^{\int_0^tf(s)\;ds}\left( \frac{1}{m}+\int_0^tr(s)\;ds\right) =  (2m^2)^{\int_0^tf(s)\;ds}\left( \frac{1}{m}+C_1\int_0^t\int_{\R^n}|\alpha_r\beta_2|\;dxds\right).
	\end{gathered}
\end{equation}
In order to complete the proof it suffices to control the part involving $\alpha_r$. A suitable decay of this term is due to the fact that the support of $\alpha_r$ shrinks with growing $m$. To take advantage of this fact additional integrability of $\alpha$ is needed. For the dissipative case this is a consequence of the energy estimates and interpolation. More precisely, we have
\begin{equation}
	\alpha \in L^q((0,T),L^p(\R^n))
\end{equation}
for some $1<q< \infty$ and $1<p<\frac{n}{n-2\gamma}$ given by the usual linear interpolation relation. For $\epsilon>0$ such that $p-\epsilon >1$, by H\"older's inequality, we get
\begin{equation}\label{16}
	\int_{\R^n}|\alpha_r\beta_2|\;dx \leq \|\alpha_r\|_{L^{p-\epsilon}}\|\beta_2\|_{L^{(p-\epsilon)'}(supp\, \alpha_r)} \leq C\|\alpha_r\|_{L^{p-\epsilon}}\|\beta_2\|_{BMO}
\end{equation}
Inequality (\ref{16}) is true for a.e. $t$ since the support of $\alpha_r$ has finite measure and, by Remark \ref{rem1},  functions with bounded mean oscillation are $p$-integrable on sets of finite measure for all $p < \infty$. From Chebyschev's inequality we obtain
\begin{equation}\label{17}
	|supp \, \alpha_r| \leq \left( \frac{\|\alpha_r\|_{L^{p}}}{m} \right)^{p} \leq \left( \frac{\|\alpha\|_{L^{p}}}{m} \right)^{p}.
\end{equation}
By H\"older's inequality 
\begin{equation}\label{18}
	\|\alpha_r\|_{L^{p-\epsilon}} = \left( \int_{supp \,\alpha_r} |\alpha_r|^{p - \epsilon}\right)^{\frac{1}{p-\epsilon}} \leq |supp \, \alpha_r|^{\frac{\epsilon}{p(p-\epsilon)}} \|\alpha_r\|_{L^{p}}
\end{equation}
Combining (\ref{17}) and (\ref{18}) we obtain
\begin{equation}
	\|\alpha_r\|_{L^{p-\epsilon}} \leq m^{-\frac{\epsilon}{p-\epsilon}} \|\alpha_r\|_{L^{p}}^{1+\frac{\epsilon}{p-\epsilon}} \leq m^{-\frac{\epsilon}{p-\epsilon}} \|\alpha\|_{L^{p}}^{1+\frac{\epsilon}{p-\epsilon}}
	\end{equation}	
Hence
\begin{equation}
	\int_{\R^n}|\alpha_r\beta_2|\;dx \leq m^{-\frac{\epsilon}{p-\epsilon}} \|\alpha\|_{L^{p}}^{1+\frac{\epsilon}{p-\epsilon}} \|\beta\|_{BMO} = m^{-\frac{\epsilon}{p-\epsilon}} \|\alpha\|_{L^{p}}^{\frac{p}{p-\epsilon}} \|\beta\|_{BMO}.
\end{equation}
Therefore
\begin{equation}\label{21}
	\int_0^t\int_{\R^n}|\alpha_r\beta_2|\;dxds \leq m^{-\frac{\epsilon}{p-\epsilon}} \|\alpha\|_{L^q((0,T),L^{p}(\R^n))} \|\beta_2\|_{L^{(q(p-\epsilon)/p)'}((0,T),BMO)}.
\end{equation}
Hence, if $\beta_2 \in  L^s((0,T),BMO)$ for some $s$ such that $q'<s<\infty$, there exists an $\epsilon>0$ such that the above conditions are satisfied and $s'=q\left(\frac{p-\epsilon}{p}\right)$. Denoting 
$$
	M = \|\beta_2\|_{L^s((0,T),BMO)},
$$
estimates (\ref{14})-(\ref{21}) give
\begin{equation}
y(t) \leq C (2m^2)^{\int_0^tf(s)\;ds}\left( \frac{1}{m}+m^{-\frac{\epsilon}{p-\epsilon}} C_1M \|\alpha\|_{L^s((0,T),L^{p}(\R^n))} \right).
\end{equation}
We choose $0\leq t_1\leq T_0$ small enough so that $2\int_0^{t_1}f(s)\;ds-1 \leq - \frac{\epsilon}{p-\epsilon}$. Then
\begin{equation}
	y(t) \leq K m^{-\frac{\epsilon}{p-\epsilon}} \;\;\;\;\; \text{ for $0 \leq t \leq t_1$,}
\end{equation}	 
where $K$ is a constant independent of $m$. Letting $m \to \infty$ we get $y(t)=0$ for $0 \leq t \leq t_1$ which implies $x(t)=0$ for $0 \leq t \leq t_1$. Therefore $\theta_1=\theta_2$ for $0 \leq t \leq t_1$. We can continue this process starting at $t=t_1$ and obtain uniqueness of $\theta$ for all $t \in [0,T]$. Uniqueness of $u$ follows from uniqueness of $\theta$.
\end{proof}

In the inviscid case $\kappa=0$, the equation does not provide any additional integrability (there is no dissipation). Therefore, in order to prove a similar theorem we have to impose this condition on our solutions. The following theorem is a corollary from Theorem \ref{thm1}.
\begin{theorem}\label{thm2}
	Let $\theta(x,t)$ be a weak solution of the equation (\ref{EQ}) such that
	$$
		\theta \in L^\infty((0,T),L^2(\R^n))\cap L^q((0,T),L^p(\R^n)) \;\;\;\;\; \text{ with $p,q \in(2,\infty)$}.
	$$
	Assume in addition $\nabla \theta =  g+b$ where $g \in L^1((0,T),L^\infty(\R^n))$ and $b \in L^s((0,T),BMO)$ for some $s>\frac{q}{q-2}$ ($s=\infty$ if $q=2$) with $supp\; b$ of finite measure for a.e. $t \in (0,T)$. Then $(\theta,u)$ is unique.
\end{theorem}

For the proof of Theorem \ref{thm1} only $L^1$-integrability of $\nabla \theta$ in time is required. Assuming that $\theta$ is not only $L^2$-integrable in space uniformly in time, but also $L^{2+\epsilon}$-integrable uniformly in time we obtain the following corollary.
\begin{corollary}
	Let $\theta(x,t)$ be a weak solution of the equation (\ref{EQ}) such that
	$$
		\theta \in L^\infty((0,T),L^2(\R^n))\cap L^q((0,T),L^p(\R^n)) \;\;\;\;\; \text{ with $p,q \in(2,\infty)$}.
	$$
	Assume in addition $\nabla \theta =  g+b$ where $g \in L^1((0,T),L^\infty(\R^n))$ and $b \in L^1((0,T),BMO)$ such that $supp\;b$ has finite measure for a.e. $t \in (0,T)$. Then $(\theta,u)$ is unique.
	
\end{corollary}
The proof of this corollary follows the same outline as the proof of Theorem 1.2 in \cite{MR08} and we refer the interested reader to this paper.

\section{Inviscid limit}

\begin{theorem}\label{thm3}
	Let $\theta_0(x,t)$ be a weak solution of (\ref{EQ}) and let $\theta_\kappa(x,t)$ be a family of weak solutions of (\ref{DEQ}) for some sufficiently small $\kappa>0$. Furthermore, assume that 
	\begin{equation}
		\begin{gathered}
			\|\theta_0(\cdot,t)\|_{L^2}+	\|\theta_\kappa(\cdot,t)\|_{L^2} \leq a_0(t) \;\;\;\;\; \text{ and $a_0 \in L^{\infty}(0,T)$}, \hfill\\
			\|\theta_0(\cdot,t)\|_{H^\gamma}+\|\theta_\kappa(\cdot,t)\|_{H^\gamma} \leq g_0(t) \;\;\;\;\; \text{ and $g_0 \in L^{2}(0,T)$}. \hfill\\
		\end{gathered}
	\end{equation} 
	Let $s>\frac{q}{q-2}$, where $q \geq 2$ ($s=\infty$ if $q=2$) is such that $\theta \in L^q((0,T),L^p(\R^n))$ for some $p \in (2,\frac{2n}{n-2\gamma})$. If  $\nabla \theta_0=g+b$ with
	\begin{equation}
		\|g(\cdot,t)\|_{L^\infty}+\|b(\cdot,t)\|_{BMO} \leq f_0(t) \;\;\;\;\; \text{ and $L^s(0,T)$}, \hfill
	\end{equation}
	where $supp\;b$ is of finite measure for a.e. $t \in(0,T)$, then 
	$$
		\sup_{0 \leq t \leq T} \|(\theta_\kappa - \theta_0)(\cdot,t)\|_{L^2} \to 0, \;\;\;\;\; \text{ as $\kappa \to 0$.}
	$$
\end{theorem}
\begin{proof}
	Let $\theta_\kappa$ and $\theta_0$ be solutions to problems (\ref{DEQ}) and (\ref{EQ}), respectively. Their difference, $\theta=\theta_\kappa-\theta_0$, satisfies
	\begin{equation}
		\theta_t - \kappa(-\Delta)^\gamma \theta_\kappa + u_\kappa\cdot \nabla \theta + (u_\kappa-u_0)\cdot\nabla \theta_0 = 0,
	\end{equation}
	where $u_\kappa$ and $u_0$ are drift velocities corresponding to $\theta_\kappa$ and $\theta_0$. Multiplying by $\theta$ and integrating in space using the divergence-free condition on drift velocities, we have 
	\begin{equation}
		\frac{1}{2}\frac{d}{dt} \int_{\R^n}|\theta|^2\;dx + \int_{\R^n} u_\kappa \cdot \nabla \theta_0 \theta\;dx + \kappa\int_{\R^n} \Lambda^\gamma \theta_\kappa \Lambda^\gamma\theta\;dx =0.
	\end{equation}
Let $\alpha = u_\kappa \theta$ and $\beta = \nabla \theta_0$. We split $\alpha =\alpha_\kappa + \alpha_r$ so that $|\alpha_\kappa| = \min(1/\kappa, |\alpha|)$. Proceeding as in the proof of Theorem \ref{thm1} we denote $f(t)=2CC_1(\|\beta_1(t)\|_{L^\infty}+\|\beta_2(t)\|_{BMO})$, $r(t)=C_1\int_{\R^n}|\alpha_r\beta_2|\;dx$ and $x(t)=C_1\|\theta\|^2_{L^2}$, and obtain the following differential inequality
\begin{equation}
	\begin{gathered}
		x'(t) \leq f(t)x(t) \left( 1+|\ln x(t)| + \ln\left(1+\frac{1}{\kappa}\right)\right) + r(t) + \kappa g^2_0(t), \\
		x(0)=0. \hfill
	\end{gathered}
\end{equation}
To find a good estimate on $x(t)$, we introduce the following equation
\begin{equation}\label{27}
	\begin{gathered}
		y'(t) =f(t)y(t) \left( 1+|\ln y(t)| + \ln\left(1+\frac{1}{\kappa}\right)\right) + r(t) + \kappa g^2_0(t), \\
		y(0)=\kappa \hfill
	\end{gathered}
\end{equation}
for some sufficiently small $\kappa$. From the Osgood existence theorem we know that there exists a unique local solution to (\ref{27}). The solution of (\ref{27}) dominates $x(t)$, i.e.
$$
	0 \leq x(t) \leq y(t) \;\;\;\;\; \text{ for $t \in [0,T_0]$,}
$$
where $T_0$ is chosen by similar rules as in the proof of Theorem \ref{thm1}. From (\ref{27}) by Gronwall's inequality we have
\begin{equation}\label{28}
	\begin{gathered}
		y(t) \leq \kappa \exp\left( \int_0^t f(s)\left[1+ |\ln y(s)| + \ln \left( 1+ \frac{1}{\kappa}\right) \right]\;ds\right)  \hfill \\
		+ \int_0^t \left( r(s)+g_0^2(s)\kappa\right)\exp\left( \int_s^t f(\tau)\left[1+ |\ln y(\tau)| + \ln \left( 1+ \frac{1}{\kappa}\right) \right]\;d\tau \right)\;ds \hfill \\
		\leq \kappa \exp \left( \left( 1+\ln\left( 1+ \frac{1}{\kappa}\right)\right) \int_0^tf(s)\;ds\right)\exp\left( \int_0^t f(s)|\ln y(s)|\;ds\right) \hfill\\
		+ \exp \left( \left( 1+\ln\left( 1+ \frac{1}{\kappa}\right)\right) \int_0^tf(s)\;ds\right)\exp\left( \int_0^t f(s)|\ln y(s)|\;ds\right)	\cdot \int_0^t \left( r(s)+g_0^2(s)\kappa\right)\;ds. \hfill
	\end{gathered}
\end{equation}
The condition $y(t) \geq \kappa$ for sufficiently small $\kappa$ gives $|\ln y(t)| \leq -\ln \kappa = \ln(1/\kappa)$. Also let $\kappa$ be small enough so that $\frac{1}{\kappa}\left(1+\frac{1}{\kappa}\right) \leq \frac{2}{\kappa^2}$, thus we can estimate the right hand-side of (\ref{28}) (up to a constant) by
\begin{equation}
	\begin{gathered}
		\kappa \left( \frac{2}{\kappa^2}\right)^{\int_0^t f(s)\;ds} + \left( \frac{2}{\kappa^2}\right)^{\int_0^tf(s)\;ds}\int_0^t(r(s)+g_0^2(s)\kappa)\;ds \hfill \\
		= \left( \frac{2}{\kappa^2}\right)^{\int_0^t f(s)\;ds} \left( \kappa + C_1\int_0^t\int_{\R^n}|\alpha_r \beta_2|\;dxds + \kappa \int_0^t g_0^2(s)\;ds\right).
	\end{gathered}
\end{equation}
Therefore we need an estimate on the part involving $\alpha_r$. As before, we use the fact that with $\kappa \to 0$, the support of $\alpha_r$ shrinks. From assumptions on $\theta_\kappa$ and $\theta_0$ we have
\begin{equation}
	\alpha \in L^\infty((0,T),L^1(\R^n)) \cap L^{1}((0,T),L^{\frac{n}{n-2\gamma}}(\R^n)),
\end{equation}
thus by interpolation
\begin{equation}
	\alpha \in L^q((0,T),L^p(\R^n)).
\end{equation}
We proceed as in the proof of Theorem \ref{thm1} and, by H\"older's inequality, have
\begin{equation}\label{32}
	\int_{\R^n}|\alpha_r\beta_2|\;dx \leq \|\alpha_r\|_{L^{p-\epsilon}}\|\beta_2\|_{L^{(p-\epsilon)'}(supp\, \alpha_r)} \leq C\|\alpha_r\|_{L^{p-\epsilon}}\|\beta_2\|_{BMO}
\end{equation}
for some sufficiently small $\epsilon >0$ (so that $p-\epsilon >1$). From Chebyschev's inequality we obtain
\begin{equation}\label{33}
	|supp \, \alpha_r| \leq \left( \kappa\|\alpha_r\|_{L^{p}} \right)^{p} \leq \left( \kappa\|\alpha\|_{L^{p}} \right)^{p}.
\end{equation}
By H\"older's inequality 
\begin{equation}\label{34}
	\|\alpha_r\|_{L^{p-\epsilon}} = \left( \int_{supp \,\alpha_r} |\alpha_r|^{p - \epsilon}\right)^{\frac{1}{p-\epsilon}} \leq |supp \, \alpha_r|^{\frac{\epsilon}{p(p-\epsilon)}} \|\alpha_r\|_{L^{p}}
\end{equation}
Combining (\ref{33}) and (\ref{34}) we obtain
\begin{equation}
	\|\alpha_r\|_{L^{p-\epsilon}} \leq \kappa^{\frac{\epsilon}{p-\epsilon}} \|\alpha_r\|_{L^{p}}^{1+\frac{\epsilon}{p-\epsilon}} \leq \kappa^{\frac{\epsilon}{p-\epsilon}} \|\alpha\|_{L^{p}}^{1+\frac{\epsilon}{p-\epsilon}}
	\end{equation}	
Hence
\begin{equation}
	\int_{\R^n}|\alpha_r\beta_2|\;dx \leq \kappa^{\frac{\epsilon}{p-\epsilon}} \|\alpha\|_{L^{p}}^{1+\frac{\epsilon}{p-\epsilon}} \|\beta_2\|_{BMO} = \kappa^{\frac{\epsilon}{p-\epsilon}} \|\alpha\|_{L^{p}}^{\frac{p}{p-\epsilon}} \|\beta_2\|_{BMO}.
\end{equation}
Therefore
\begin{equation}\label{37}
	\int_0^t\int_{\R^n}|\alpha_r\beta_2|\;dxds \leq \kappa^{\frac{\epsilon}{p-\epsilon}} \|\alpha\|_{L^q((0,T),L^{p}(\R^n))} \|\beta_2\|_{L^{(q(p-\epsilon)/p)'}((0,T),BMO)}.
\end{equation}
Thus, if we have $\beta_2 \in L^s((0,T),BMO)$ for some $s$ such that $q' < s < \infty$, there exists an $\epsilon>0$ such that above conditions are satisfied and $s'=q\left(\frac{p-\epsilon}{p}\right)$. We denote
$$
	M = \|\beta_2\|_{L^s((0,T),BMO)}.
$$
Then, estimates (\ref{28})-(\ref{37}) give
\begin{equation}
y(t) \leq C \left( \kappa^{1-2\int_0^tf(s)\;ds}+\kappa^{\frac{\epsilon}{p-\epsilon}} C_1M \|\alpha\|_{L^q((0,T),L^{p}(\R^n))} \right).
\end{equation}
We choose $0\leq t_1\leq T_0$ small enough so that $1-2\int_0^{t_1}f(s)\;ds \leq \frac{\epsilon}{p-\epsilon}$. We get
\begin{equation}
	y(t) \leq K \kappa^{\frac{\epsilon}{p-\epsilon}} \;\;\;\;\; \text{ for $0 \leq t \leq t_1$,}
\end{equation}	 
where $K$ is a constant independent of $\kappa$. We can continue the process on an interval starting at $t=t_1$ but with the initial condition in (\ref{27}) changed to $y(t_1)=K\kappa^{\frac{\epsilon}{p-\epsilon}}$. Due to integrability of $f(t)$, iterating this procedure, we eventually cover the whole interval $[0,T]$. Letting $\kappa \to 0$ we obtain
\begin{equation}
	\sup_{0\leq t \leq T}\|\theta_\kappa - \theta_0\|_{L^2} \to 0, \;\;\;\;\; \text{ as $\kappa \to 0$}.
\end{equation}
\end{proof}
\begin{remark}
	From the proof of Theorem \ref{thm3} follows that if we assume 
	$$
		\sup_{0 \leq t \leq T}	|f_0(t)+g_0^2(t)| \leq N,
	$$
	then we obtain the following explicit rate of convergence
	$$
		\sup_{0 \leq t \leq T}\|\theta_\kappa - \theta_0\|_{L^2} \leq C\kappa^{e^{-2NT}}.
	$$
	The proof follows as in \cite{MR08}.
\end{remark}

\appendix

\section{Existence of weak solutions to (\ref{DEQ})}

We sketch the proof of existence of global Leray-Hopf weak solutions of equations (\ref{DEQ}). We follow the general outline used to construct weak solutions of the Navier-Stokes equations (cf. \cite{Tem}). Let $\phi \in C^\infty_0(\R^n)$ be positive, with $\int_{\R^n}\phi \;dx=1$. Then $\phi_\epsilon(x) = \epsilon^{-n}\phi(x/\epsilon)$, for $\epsilon >0$, is a standard family of mollifiers. We consider the approximating system
\begin{equation}\label{app1}
	\begin{gathered}
		\partial_t \theta^\epsilon + (u^\epsilon \cdot \nabla)\theta^\epsilon + \kappa(-\Delta)^\gamma\theta^\epsilon = \epsilon \Delta^2\theta^\epsilon \hfill
	\end{gathered}
\end{equation}
with initial data $\theta_0^\epsilon = \phi_\epsilon \ast \theta_0$. The drift velocity $u^\epsilon$ is obtained from $\theta^\epsilon$ via a Calder\'on-Zygmund singular integral operator, defined as a Fourier multiplier with a zero-order symbol. Note that $\|\theta_0^\epsilon\|_{L^2} \leq \|\theta_0\|_{L^2}$ for any $\epsilon >0$. 

Let $s > \frac{n}{2} + 1$ and fix $\epsilon >0$. Notice that $(-\Delta)^{s/2}\theta_0 \in L^2(\R^n)$. Since $\epsilon \Delta^2 $ provides enough dissipation, from standard energy estimates it follows that 
$$
	\sup_{t \in [0,T]}\|\theta^\epsilon(t)\|_{H^s} \leq C_{\epsilon,n,T,\phi,\|\theta_0\|_{L^2}},
$$ 
where $C_{\epsilon,n,T,\phi,\|\theta_0\|_{L^2}}$ is a strictly positive constant which is finite for any $0<T < \infty$. This a priori estimate and a classical Galerkin approach ensure the global existence of a strong $H^s$ solution to (\ref{app1}). Moreover, we have the energy estimate
\begin{equation}\label{app2}
	\|\theta^\epsilon(t)\|_{L^2} + \kappa\int_0^t\|\theta^\epsilon(s)\|^2_{H^\gamma}\;ds \leq \|\theta_0\|^2_{L^2},
\end{equation}
which is uniform in $\epsilon >0$. Thus
\begin{equation}\label{app3}
	\theta^\epsilon \in C([0,T],L^2(\R^n)) \cap L^2((0,T),H^\gamma(\R^n)).
\end{equation}
This guarantees that, up to a subsequence, $\theta^\epsilon$ converges weakly-$\ast$ and weakly to some function $\theta$ in the spaces $L^\infty((0,T),L^2(\R^n))$ and $L^2((0,T),H^\gamma(\R^n))$, respectively. This however is not enough to conclude that $\theta$ solves the equation (\ref{DEQ}) or (\ref{EQ}). The classical condition that allows one to conclude that $\theta$ is indeed a solution is to establish uniform control over $\partial_t \theta^\epsilon$ in some negative Sobolev space. This is accomplished using the energy inequality (\ref{app2}). From the Aubin-Lions compactness lemma we obtain that 
$$
	\theta^\epsilon \to \theta \text{ strongly in $L^2((0,T),L^2_{loc}(\R^n))$.}
$$
The strong convergence is necessary in order to pass to the limit in the nonlinear term in the weak formulation of (\ref{DEQ}). This is done in a standard way using the fact that Calder\'on-Zygmund operators map $L^2 \to L^2$ strongly.
	
	\bigskip
	
\noindent\footnotesize{\bf{Acknowledgements.}} We wish to express gratitude to Susan Friedlander for fruitful discussions and comments on early versions of this draft. 

\end{document}